\documentclass[12pt]{amsart}

\usepackage{amssymb}

\usepackage{enumerate}

\usepackage{graphicx}

\makeatletter
\@namedef{subjclassname@2010}{%
  \textup{2010} Mathematics Subject Classification}
\makeatother

\newtheorem{theorem}{ Main Theorem}[section]

\theoremstyle{definition}

\numberwithin{equation}{section}

\begin{document}


\baselineskip=17pt


\title[On the Diophantine equation $ \sum_{i=1}^n a_ix_{i} ^4= \sum_{j=1}^na_j y_{j}^4 $ ]{On the Diophantine equation  $ \sum_{i=1}^n a_ix_{i} ^4= \sum_{j=1}^na_j y_{j}^4 $}

\author[F. Izadi]{Farzali Izadi}
\address{Farzali Izadi \\
Department of Mathematics \\ Faculty of Science \\ Urmia University \\ Urmia 165-57153, Iran}
\email{f.izadi@urmia.ac.ir}

\author[M. Baghalaghdam]{Mehdi Baghalaghdam}
\address{Mehdi Baghalaghdam \\
Department of Mathematics\\ Faculty of Science \\ Azarbaijan Shahid Madani University\\Tabriz 53751-71379, Iran}
\email{mehdi.baghalaghdam@azaruniv.edu}

\date{}

\begin{abstract}
In this paper, by using the elliptic curves theory, we study the fourth power Diophantine equation ${ \sum_{i=1}^n a_ix_{i} ^4= \sum_{j=1}^na_j y_{j}^4 }$, where $a_i$ and $n\geq3$ are fixed arbitrary integers. We solve the equation for some values of $a_i$ and $n=3,4$, and find nontrivial solutions for each case in natural numbers. By our method, we may find infinitely many nontrivial solutions for the above Diophantine equation and show, among the other things, that how some numbers can be written as sums of three, four, or more biquadrates in two different ways. While our method can be used for solving  the equation  for every $a_i$ and $n\geq 3$, this paper will  be
restricted  to the examples where $n=3,4$. In the end, we explain how to solve it in general cases without giving concrete examples.
\end{abstract}



\maketitle

\section{Introduction}
Number theory is a vast and fascinating field of mathematics, sometimes called ''higher arithmetic'', consisting
of the study of the properties of whole numbers. In mathematics, a Diophantine equations  is a polynomial equation, usually in two or more unknowns, such that only the integer solutions are studied. The word Diophantine refers to the Hellenistic mathematician of the 3rd century, Diophantus of Alexanderia, who made a study of such equations. While individual equation present  a kind of puzzle and have been considered throughout the history, the formulation of general theories of Diophantine equations (beyond the theory of quadratic forms) was an achievement of the twentieth century. As a history work, Euler conjectured in $1969$ that the Diophantine equation $A^4+B^4+C^4=D^4$, or more generally $A_1^N+A_2^N+ \cdots +A_{N-1}^N=A_N^N$, ($N\geq4$), has no solution in positive integers (see \cite{L.D}). Nearly two centuries later, a computer search (see \cite{L.L}) found the first counterexample to the general conjecture (for $N=5$):
$27^5+84^5+110^5+133^5=144^5$.

 In $1986$, Noam Elkies found a method to construct an infinite series of counterexamples for the $K=4$ case (see \cite{N.E}). His smallest counterexample was:

 $2682440^4+15365639^4+18796760^4=20615673^4$.

 In this paper, we are interested in the study of the Diophantine equation:

${ \sum_{i=1}^n a_ix_{i} ^4= \sum_{j=1}^na_j y_{j}^4 }$,

where $a_i$, $n\geq 3$ are  fixed arbitrary integers.

Now we are going to solve the equation for the cases $n=3,4$, and in the end we explain how  to solve it in general cases ($n\geq 5$).
\\
Our main results are the followings:

\begin{theorem} Consider the Diophantine equation :
\\
$ax^4+by^4+cz^4=au^4+bv^4+cw^4$,
where all the coefficients are fixed arbitrary integers. Let $Y^2=X^3+fX^2+gX+h$, be an elliptic curve in which the coefficients of $f$, $g$, and $h$ are all functions
of $a$, $b$, $c$, and two other rational parameters $A$  and $B$, yet to be found later. If the elliptic curve has positive rank, depending on the values of $A$ and $B$, the Diophantine equation has infinitely many  integer solutions.
\end{theorem}

Proof. Let: $x=m+p$, $y=m-q$, $z=m-s$, $u=m-p$, $v=m+q$ and $w=m+s$, where all variables are rational numbers. By substituting these variables in the above Diophantine equation, and after some simplifications, we get:
\\
$m^2(ap-bq-cs)=-ap^3+bq^3+cs^3$.
\\

We may assume that $ap-bq-cs=1$ and $m^2=-ap^3+bq^3+cs^3$. Also let $q=As+B$, where $A$, $B$ $\in \mathbb{Q}$. By plugging
$q=As+B$  and $p=\frac{b}{a}q+\frac{c}{a}s+\frac{1}{a}$,
 into the equation $m^2=-ap^3+bq^3+cs^3$,  and after some simplifications, we obtain the elliptic curve:
\\

$m^2=(-a(\frac{bA+c}{a})^3+bA^3+c)s^3+(-3a(\frac{bA+c}{a})^2.(\frac{bB+1}{a})+3A^2Bb)s^2
\\
+
(-3a(\frac{bA+c}{a}).(\frac{bB+1}{a})^2+3AB^2b)s+(-a(\frac{bB+1}{a})^3+bB^3)$.
\\

 Multiplying the both sides of this elliptic curve by

  $(-a(\frac{bA+c}{a})^3+bA^3+c)^2$ and letting
\begin{equation} \label{2}
Y=(-a(\frac{bA+c}{a})^3+bA^3+c)m,
\end{equation}
and
\begin{equation} \label{3}
 X=(-a(\frac{bA+c}{a})^3+bA^3+c)s,
\end{equation}

we get the new elliptic curve $Y^2=X^3+fX^2+gX+h$, where\\
\begin{equation*}
f=(-3a(\frac{bA+c}{a})^2.(\frac{bB+1}{a})+3A^2Bb),
\end{equation*}
 \\
 \begin{equation*}
g=(-3a(\frac{bA+c}{a}).(\frac{bB+1}{a})^2+3AB^2b).(-a(\frac{bA+c}{a})^3+bA^3+c),
\end{equation*}
\\
\begin{equation*}
h=(-a(\frac{bB+1}{a})^3+bB^3).(-a(\frac{bA+c}{a})^3+bA^3+c)^2.
\end{equation*}
 \\

If the above elliptic curve has positive rank (for every $a$, $b$, $c$, this is done  by choosing appropriate values for $A$ and $B$.), by calculating $m$, $s$, $q$, $p$, $m\pm p$, $m\pm q$, $m\pm s$, from the relations \eqref{2}, \eqref{3}, $q=As+B$  and $p=\frac{b}{a}q+\frac{c}{a}s+\frac{1}{a}$ , and after some simplifications  and canceling  the  denominators of  $m\pm p$, $m\pm q$, $m\pm s$, we may obtain infinitely many integer solutions for the above Diophantine equation. Now the proof is complete.
\\

Now we are going to solve some couple of examples:
\\
 Example 1. We wish to solve the Diophantine equation:
 \\
$x^4+2y^4+3z^4=u^4+2v^4+3w^4$,
\\
We may assume that $A=4$, $B=0$. (we may choose the other appropriate values for $A$ and $B$ so that the corresponding elliptic curve has positive rank).

Then we get the elliptic curve:

 $Y^2=X^3-363X^2+39600X-1200^2$.

The rank of this elliptic curve is $1$ and its generator is the point

 $P=(X,Y)=(\frac{3625}{16},\frac{46525}{64})$. Because of this, the above elliptic curve has infinitely many rational points and we may obtain infinitely many solutions for the Diophantine equation too. Since $X=\frac{3625}{16}$ and $Y=\frac{46525}{64}$, by calculating $m$, $s$, $q$, $p$, $m\pm p$, $m\pm q$, $m\pm s$, from the relations \eqref{2}, \eqref{3} and $q=As+B$, $p=\frac{b}{a}q+\frac{c}{a}s+\frac{1}{a}$, and after some simplifications and canceling the denominators of $m\pm p$, $m\pm q$, $m\pm s$, we get the identity:
 \\

$5169^4+2.(459)^4+3.(1281)^4=1447^4+2.(4181)^4+3.(2441)^4$.
\\

 By choosing the other points on the elliptic curve such as $2P$, $3P$, ..., we obtain infinitely many solutions for the  Diophantine equation.
\\

Example 2. Let $a=1$, $b=1$, $c=61$, then by letting $A=2$, $B=0$, we get the elliptic curve
\\
$Y^2=X^3-11907X^2+47245842X-249978^2$.
\\
The rank of this elliptic curve is $1$ and its generator is the point

 $P=(X,Y)=(\frac{2613213887380271422}{612348332222929},\frac{-35386313782867169078293498}{15152971591283964136217})$.
\\

 By calculating the above values and after some simplifications, we obtain a solution for the Diophantine equation:  $x^4+y^4+61z^4=u^4+v^4+61w^4$ as
 \\

$183488684443575775594469^4+120584031079948181257985^4
\\
+61.(73244546207202190584444)^4=
235298807112488175416275^4
\\
+68773908411035781436179^4+61.(21434423538289790762638)^4$.
\\

By choosing the other points on the elliptic curve such as $2P$, $3P$, ..., we obtain infinitely many solutions for the Diophantine equation.
\\

Example 3. Sums of three biquadrates in two different ways:
\\
In this case we have $a=b=c=1$. Letting $A=-10$, $B=0$, we get the elliptic curve
$Y^2=X^3-243X^2-7290X-270^2$.
\\
The rank of this elliptic curve is $1$ and its generator is the point

 $P=(X,Y)=(450,6210)$. By calculating the above values and after some simplifications, we obtain a solution for the Diophantine equation:\\ $x^4+y^4+z^4=u^4+v^4+w^4$ as

$19^4+74^4+117^4=21^4+64^4+119^4$.
\\

Also we have: $2P=(X',Y')=(\frac{606357}{2116},\frac{-115780401}{97336})$ and
\\

 $3P=
(X'',Y'')=(\frac{255306774610}{164070481},\frac{-118288360159623370}{2101578791129})$.
\\

By using these two new points, we obtain the two other solutions for the Diophantine equation respectively:
\\

$17948013^4+43856069^4+9765331^4=43676991^4+18127091^4+15963647^4$,
\\

$8828891360220313^4+15099060491941827^4+11501813568364388^4=
\\
14828780671704361^4+8558611539982847^4+12155858463560286^4$.
\\

Choosing the other points such as $4P$, $5P$, ..., give rise to infinitely many solutions for the  Diophantine equation.
\\

\begin{theorem} Consider the Diophantine equation:
\\
$ax^4+by^4+cz^4+dt^4=au^4+bv^4+cw^4+dh^4$,

where all the coefficients are fixed arbitrary integers.

 Let $Y^2=X^3+fX^2+gX+h$, be an elliptic curve in which the coefficients of $f$, $g$, and $h$ are all functions
of $a$, $b$, $c$, $d$, and four other rational parameters $A$, $B$, $D$, $F$, yet to be found later. If the elliptic curve has positive rank, depending on the values of $A$, $B$, $D$, $F$, the Diophantine equation has infinitely many integer solutions.
\end{theorem}

Proof. Let $x=m+p$, $y=m-q$, $z=m-s$, $t=m+r$, $u=m-p$, $v=m+q$, $w=m+s$ and $h=m-r$, where all variables are rational numbers. By substituting these variables in the above Diophantine equation, and after some simplifications, we get:
\\

$m^2(ap-bq-cs+dr)=-ap^3+bq^3+cs^3-dr^3$.
\\

we may assume that $ap-bq-cs+dr=1$ and $m^2=-ap^3+bq^3+cs^3-dr^3$. Also let $q=As+B$ and $r=Ds+F$, where $A$, $B$, $D$, $F$ $\in \mathbb{Q}$. By plugging $q=As+B$, $r=Ds+F$ and  $p=\frac{b}{a}q+\frac{c}{a}s-\frac{d}{a}r+\frac{1}{a}$ into the equation $m^2=-ap^3+bq^3+cs^3-dr^3$, and after some simplifications, we obtain the elliptic curve:
 \\

$m^2=(-a(\frac{bA+c-dD}{a})^3+bA^3+c-dD^3)s^3+(-3a(\frac{bA+c-dD}{a})^2.(\frac{bB-dF+1}{a})+3A^2Bb-3dD^2F)s^2+
(-3a(\frac{bA+c-dD}{a}).(\frac{bB-dF+1}{a})^2+3AB^2b-3DF^2d)s+(-a(\frac{bB-dF+1}{a})^3+bB^3-dF^3)$.
\\

Multiplying the both sides of this elliptic curve by

$(-a(\frac{bA+c-dD}{a})^3+bA^3+c-dD^3)^2$ and letting
 \\
 \begin{equation}\label{4}
Y=(-a(\frac{bA+c-dD}{a})^3+bA^3+c-dD^3)m,
\end{equation}
 and
 \begin{equation}\label{5}
X=(-a(\frac{bA+c-dD}{a})^3+bA^3+c-dD^3)s,
\end{equation}
 \\

we get  the elliptic curve
$Y^2=X^3+fX^2+gX+h$, where
\\
\begin{equation*}
f=(-3a(\frac{bA+c-dD}{a})^2.(\frac{bB-dF+1}{a})
+3A^2Bb-3dD^2F),
\end{equation*}

\begin{equation*}
g=(-3a(\frac{bA+c-dD}{a}).(\frac{bB-dF+1}{a})^2+3AB^2b-3DF^2d)
\end{equation*}
\begin{equation*}
(-a(\frac{bA+c-dD}{a})^3+bA^3+c-dD^3),
\end{equation*}

\begin{equation*}
h=(-a(\frac{bB-dF+1}{a})^3+bB^3-dF^3).(-a(\frac{bA+c-dD}{a})^3+bA^3+c-dD^3)^2.
\end{equation*}
\\

If the above elliptic curve has positive rank, then by calculating $m$, $s$, $q$, $p$, $r$, $m\pm p$, $m\pm q$, $m\pm s$, $m\pm r$, from the relations \eqref{4}, \eqref{5} and $q=As+B$, $r=Ds+F$, $p=\frac{b}{a}q+\frac{c}{a}s-\frac{d}{a}r+\frac{1}{a}$, after  some simplifications  and canceling  the  denominators of  $m\pm p$, $m\pm q$, $m\pm s$, $m\pm r$, we may obtain infinitely many integer solutions for the  Diophantine equation. The proof is complete.
\\

Now we are going to solve some couple of examples:

Example 1. Let us to solve the Diophantine equation:
\\
 $x^4+y^4+z^4+19t^4=u^4+v^4+w^4+19h^4$,
 \\
Suppose that $A=2$, $B=0$, $D=4$, $F=0$.

Then we get the elliptic curve:
\\
 $Y^2=X^3-15987X^2+84930390X-387810^2$.
 \\
The rank of this elliptic curve is $1$ and its generator is the point

 $P=(X,Y)=(\frac{8832851584572306}{887637201025},\frac{-260518741182457285866354}{836282950759698625})$.
\\

Because of this,  the  above elliptic curve has infinitely many rational points and this gives rise infinitely many solutions for the Diophantine equation too.  By  using the point $P =(X,Y)$ and calculating $m$, $s$, $q$, $p$, $r$, $m\pm p$, $m\pm q$, $m\pm s$, $m\pm r$, from the above relations  and after some simplifications, we get the identity:
 \\

$2923081816382045453549^4+1490120403735220625479^4
\\
+1445379398594646027434^4+
19.(1221674372891773037209)^4=
\\
121805029473902594771^4+1311156383172922233299^4
\\
+
1355897388313496831344^4+19.(1579602414016369821569)^4$.
\\

By choosing the other points on the elliptic curve  such as $2P$, $3P$, ..., we obtain infinitely many solutions for the  Diophantine equation.
\\

Example 2. Let $a=1$, $b=c=d=1000$. By taking $A=2$, $B=0$, $D=3$, $F=0$, we get the elliptic curve

$Y^2=X^3-18000^2$.

The rank of this elliptic curve is $1$ and its generator is the point

$P=(X,Y)=(1000,26000)$. By calculating the above values and after some simplification, we obtain a solution for the Diophantine equation as

$x^4+1000y^4+1000z^4+1000t^4=u^4+1000v^4+1000w^4+1000h^4$,
\\

$8^4+1000.(24)^4+1000.(25)^4+1000.(29)^4=44^4+1000.(23)^4+1000.(27)^4+1000.(28)^4.$
\\

For the above solution it is interesting to see that   

 $24+25+29=23+27+28$.
Again by choosing the other points such as $2P$, $3P$, ..., one can obtain infinitely many solutions for the  Diophantine equation.
 \\

Example 3. Sums of four biquadrates in two different ways.
\\
In this case for the coefficients of the diophantine equation we have $a=b=c=d=1$. By letting $A=3$, $B=0$, $D=7$, $F=0$, we get the elliptic curve:
\\
$Y^2=X^3-27X^2-2592X-288^2$.
\\
The rank of this elliptic curve is $1$ and its generator is the point

 $P=(X,Y)=(328,5608)$. By calculating the above values and after some simplifications, we obtain a solution for the Diophantine equation as
\\
$x^4+y^4+z^4+t^4=u^4+v^4+w^4+h^4$,

$271^4+289^4+330^4+494^4=207^4+371^4+412^4+430^4$.
\\

Also we have: $2P=(X',Y')=(\frac{48232180}{491401},\frac{203244176836}{344472101})$ and
\\

 $3P=
(X'',Y'')=(\frac{61771836632004160}{797318963764569},
\frac{-3119953836085429109330528}{22513765616829041228253})$.
\\

By using these two new points, we obtain the other solutions for the Diophantine equation respectively:
\\

$325492151^4+12726487787^4+21179177332^4+54989935512^4=
\\
50485552058^4+38084556422^4+
29631866877^4+4178891303^4$,
\\

$463645068132430760337286^4+261021177580969389283009^4+
\\
152006097445436236205389^4+
284054223096696376105091^4=
\\
268647953377091441004128^4+66024062825630069949851^4+
\\
42991017309903083127769^4+479051337852035695438249^4$.
\\

By choosing the other points on the elliptic curve  such as  $4P$, $5P$, ..., we obtain  infinitely many  solutions for the  Diophantine equation.
 \\
In the end, we prove the general result without without discussing any concrete examples.
\begin{theorem} Consider the Diophantine equation:
\\
$ \sum_{i=1}^n a_ix_{i} ^4= \sum_{j=1}^na_j y_{j}^4 $,

where all the coefficients are fixed arbitrary integers and $n\geq 5$ .

 Let $Y^2=X^3+fX^2+gX+h$, be an elliptic curve in which the coefficients of $f$, $g$, and $h$ are all functions
of $a_i$, and the other rational parameters $A_i$, $B_i$, yet to be found later. If the elliptic curve has positive rank, depending on the values of $A_i$, $B_i$,  the Diophantine equation has infinitely many  integer solutions.
\end{theorem}

 Proof. Let: $x_i= m+p_i$ and $y_i=m-p_i$, (For $ 1\leq i \leq n$), where $m$ and $p_i$  $\in \mathbb{Q}$. After substituting these variables in the Diophantine equation, we obtain:
\\

 $m^2(a_1p_1+ \cdots +a_np_n)=-a_1p_1^3- \cdots -a_np_n^3$.
\\

We may assume that $(a_1p_1+ \cdots +a_np_n)=1$ and  $m^2=-a_1p_1^3- \cdots -a_np_n^3$.
 Also we let $p_i=A_ip_n+B_i$ (for $ 2\leq i \leq n-1$).
\\

By plugging  $p_1=-\frac{a_2}{a_1}p_2- \cdots-\frac{a_{n-1}}{a_1}p_{n-1} -\frac{a_n}{a_1}p_n+\frac{1}{a_1}=Gp_n+H$, where
\\

 $G=\frac{-a_2A_2-a_3A_3- \cdots -a_{n-1}A_{n-1}-a_n}{a_1}$,
 \\
 and

 $H=\frac{-a_2B_2-a_3B_3- \cdots -a_{n-1}B_{n-1}+1}{a_1}$,
\\

and  $p_i=A_ip_n+B_i$ (for $ 2\leq i \leq n-1$), into the equation

  $m^2=-a_1p_1^3- \cdots -a_np_n^3$, and after some simplifications, we obtain the elliptic curve
\\

$m^2=L_1p_n^3+L_2p_n^2+L_3p_n+L4$, where
\\

$L_1=(-a_1G^3-a_2A_2^3-a_3A_3^3- \cdots - a_{n-1}A_{n-1}^3-a_n)$,
\\

$L_2=(-3a_1G^2H-3a_2A_2^2B_2-3a_3A_3^2B_3- \cdots -3a_{n-1}A_{n-1}^2B_{n-1})$,
\\

$L_3=(-3a_1GH^2-3a_2A_2B_2^2-3a_3A_3B_3^2- \cdots -3a_{n-1}A_{n-1}B_{n-1}^2)$,
\\

$L_4=(-a_1H^3-a_2B_2^3-a_3B_3^3- \cdots- a_{n-1}B_{n-1}^3)$.
\\

Multiplying the both sides of this elliptic curve by

  $L_1^2=(-a_1G^3-a_2A_2^3-a_3A_3^3- \cdots - a_{n-1}A_{n-1}^3-a_n)^2$ and letting
\begin{equation} \label{6}
Y=(-a_1G^3-a_2A_2^3-a_3A_3^3- \cdots - a_{n-1}A_{n-1}^3-a_n)m,
\end{equation}
and
\begin{equation} \label{7}
 X=(-a_1G^3-a_2A_2^3-a_3A_3^3- \cdots - a_{n-1}A_{n-1}^3-a_n)p_n,
\end{equation}
\\

we get the elliptic curve $Y^2=X^3+fX^2+gX+h$, where\\

$f=(-3a_1G^2H-3a_2A_2^2B_2-3a_3A_3^2B_3- \cdots -3a_{n-1}A_{n-1}^2B_{n-1})$,
\\

$g=(-3a_1GH^2-3a_2A_2B_2^2-3a_3A_3B_3^2- \cdots -3a_{n-1}A_{n-1}B_{n-1}^2)$

$(-a_1G^3-a_2A_2^3-a_3A_3^3- \cdots - a_{n-1}A_{n-1}^3-a_n),$
\\

$h=(-a_1H^3-a_2B_2^3-a_3B_3^3- \cdots- a_{n-1}B_{n-1}^3)$

$(-a_1G^3-a_2A_2^3-a_3A_3^3- \cdots - a_{n-1}A_{n-1}^3-a_n)^2.$
\\

Then if the above elliptic curve has positive rank (this is done  by choosing appropriate values for $A_i$ and $B_i$.), by calculating $m$, $p_i$, $m\pm p_i$, from the relations \eqref{6}, \eqref{7} and  $p_i=A_ip_n+B_i$ (for $ 2\leq i \leq n-1$), $p_1=Gp_n+H$, and after some simplifications  and canceling  the  denominators of  $m\pm p_i$, we may obtain infinitely many integer solutions for the  Diophantine equation. The proof is complete.
\\

Finally we mention that each point on the elliptic curve can be represented in the form of $(\frac{r}{s^2},\frac{t}{s^3})$, where $r$, $s$, $t$ $\in \mathbb{Z}$.\\
Then for every generator $P$ of the curve we have
$nP=(\frac{r_n}{s_n^2},\frac{t_n}{s_n^3})$. This gives rise to a parametric solution for each case of the Diophantine equations. Also by choosing  the other appropriate values of $A_i$ and $B_i$ and getting the new elliptic curve of rank $\geq 1$( and repeating the above process), we may obtain infinitely many nontrivial parametric solutions for each case of Diophantine equations.


\begin{thebibliography}{HD}
\bibitem{L.D} L. E. DICKSON, History  of the  Theory of Numbers, Vol. II: Diophantine Analysis, G. E. Stechertk  Co., New York,  $1934$.
\bibitem{N.E} N. ELKIES, $ "On  A4 + B4 + C4 = D4"$. Mathematics of Computation. $1988$. $51 (184): 825–835.$
\bibitem{L.L} L. J. LANDER and T. R. PARKIN, "Counterexamples to Euler's conjecture on sums of like
powers, "Bull. Amer. Math. Soc, Vol. $72, 1966, p. 1079$.










\normalsize
\baselineskip=17pt



\bibitem[]{}

\bibitem[]{}

\bibitem[]{}

\end{thebibliography}
\end{document}